\documentclass[reqno]{amsart}
\usepackage{hyperref}
\newtheorem{theorem}{Theorem}[section]
\newtheorem{proposition}{Proposition}[section]

\newtheorem{corol}[theorem]{Corollary}

\theoremstyle{remark}

\numberwithin{equation}{section}

\begin{document}

\title[MALLIAVIN-STEIN METHOD FOR MULTI-DIMENSIONAL POISSON U-STATISTICS]{MALLIAVIN-STEIN METHOD FOR MULTI-DIMENSIONAL U-STATISTICS OF POISSON POINT PROCESSES}
\author{Nguyen Tuan Minh}
\address{Faculty of Applied Mathematics and Computer Science, Belarusian State University, 4, F. Skorina Ave., 220050 Minsk, Belarus}
\email{minhnguyen@yandex.ru}

\begin{abstract}In this paper, we give an upper bound for a probabilistic distance between a Gaussian vector and a vector of U-statistics of Poisson point processes by applying Malliavin-Stein inequality on the Poisson space.
\end{abstract}

\maketitle
\tableofcontents
\section{Introduction}
The theory of Malliavin calculus on the Poisson space was firstly studied by Nualart and Vives in their excellent paper of Strasbourg's seminars \cite{Nualart1990}. For some important contributions and applications for Poisson point processes, we refer to \cite{Houdre1995, Wu2000, Peccati2011, Last2011}. The combination of Stein's method and Malliavin calculus on the Poisson space related to the normal approximations of Poisson functionals has been considered by Peccati, Sol\'e, Taqqu, Utzet and Zheng in their recent papers \cite{Peccati2010a, Peccati2010b}.\\

The basic theory of U-statistics was introduced by Hoeffding \cite{Hoeffding1948} as a class of statistics that is especially important in estimation theory. Further applications are widely regarded in theory of random graphs, spatial statistics, theory of communication and stochastic geometry, see e.g. \cite{Lee1990, Koroljuk1994,Borovskikh1996}.\\

Recently, the idea of central limit theorems for U-statistics of Poisson point processes by using Malliavin calculus and Stein's method has been given by Reitzner and Schulte in \cite{Reitzner2011}. Our main work in the current paper is to extend their results to vectors of U-statistics by applying the muti-dimensional Malliavin-Stein inequality, that was proved by Peccati and Zheng in \cite{Peccati2010b}. Some preliminaries of Malliavin calculus on the Poisson space and U-statistic will be introduced in Section 2. An upper bound for a probabilistic distance between a Gaussian vector and a square integrable random variable with finite Wiener-It\^o chaos expansions will be shown in Section 3 and its application for multi-dimensional U-statistics of Poisson point processes will be given in Section 4.

\section{Malliavin Calculus on the Poisson space}
Let $(E, \mathcal A, \mu)$ be some measure space with $\sigma$-finite measure $\mu$. The Poisson point process or Poisson measure with intensity measure $\mu$ is a family of random variables $\{{N}(A)\}_{A\in\mathcal{A}}$ defined on some probability space $(\Omega, \mathcal F, \mathrm{P})$ such that
\begin{enumerate}
\item[1.] $ \forall A\in\mathcal{A}$, $N(A)$ is a Poisson random variable with rate $\mu(A)$.
\item[2.] If sets $A_1,A_2,\ldots,A_n\in\mathcal{A}$ don't intersect then the corresponding random variables ${N}(.)$ from i) are mutually independent.
\item[3.] $\forall\omega\in\Omega$, ${N}(.,\omega)$ is a measure on $(E, \mathcal A)$.
\end{enumerate}
Note that, for some $\sigma$-finite measure space $(E, \mathcal A, \mu)$, we can always set
$$\Omega=\{\omega =\sum_{j=1}^n\delta_{z_j}, \ n\in\mathbb{N}\cup\{\infty\}, z_j\in E  \},$$
where $\delta_z$ denotes the Dirac mass at $z$, and for each $A\in \mathcal{A}$, we give the mapping ${N}$ such that
$$\omega\mapsto {N}(A,\omega)=\omega(A).$$
Moreover, the $\sigma$-field $\mathcal{F}$ is supposed to be the $\text{P}$-completion of the $\sigma$-field generated by ${N}$.\\

Let $\mathcal{M}(E)$ denote the space of all integer-valued $\sigma$-finite measures on $E$, which can be equipped with the smallest $\sigma$-algebra $\Sigma$ such that for each $A\in\mathcal{A}$ then the mapping $\eta\in\mathcal{M}(E)\mapsto \eta(A)$  is measurable. We can give on $(\mathcal{M}(E),\Sigma)$ the probability measure $P_{N}$  induced by the Poisson measure $N$  and denote $L^p(P_{N})$ as the set of all measurable functions $F:\mathcal{M}(E)\to \overline{\mathbb R}$ such that $\mathbf{E}[|F|^p]<\infty$, where the expectation takes w.r.t. the probability measure $P_{N}$.\\

Let $L^p(\mu^n)$  be the space of all measurable functions $f: E^n \to\overline{\mathbb{R}}$  such that
$$\|f\|=\int\limits_{E^k}|f(z_1,\hdots,z_n)|^p \, \mu^n(dz_1, \dots ,dz_n)<\infty.$$
Note that $L^2(\mu^n)$ becomes a Hilbert space when we define on it the scalar product
$$ \langle f, g\rangle =\int\limits_{E^k}f(z_1,\hdots,z_n) g(z_1,\hdots,z_n)\, \mu^n(dz_1, \dots ,dz_n).$$
We denote $L^p_{\text{sym}}(\mu^n)$ as the subset of symmetric functions in $L^p(\mu^n)$, in the sense that the functions are invariant under all permutations of the arguments. Now, for each $A\in\mathcal{A}$, we define the random variable $\widehat{N}(A)=N(A)-\mu(A)$, which is also known as a compensated Poisson measure.
For each symmetric function $f\in L^2_{\text{sym}}(\mu^n)$, one can define the multiple Wiener-It\^o integral $I_n(f)$ w.r.t the compensated Poisson measure $\widehat{N}$ denoted by
\begin{equation}\label{integral}I_n(f)=\int_{E^n}f(z_1,\hdots, z_n) \widehat{N}^n(dz_1, \dots ,dz_n).\end{equation}
At first, let $\mathcal{S}_n$ be the class of simple functions $f$, which takes the form
\begin{equation}\label{simple}f(z_1,z_2,....,z_n)=\sum_{k=1}^m\lambda_k \mathbf{1}_{A_{1}^{(k)}\times\ldots \times A_{n}^{(k)}}(z_1,z_2,....,z_n),\end{equation}
where $A_{i}^{(k)}\in\mathcal{A}, \lambda_k\in \mathbb{R}$ and the sets $A_{1}^{(k)}\times\ldots \times A_{n}^{(k)}$ are pairwise disjoint such that $f$ is symmetric and vanishes on diagonals, that means $f(z_1, \ldots , z_n) = 0$ if $z_i = z_j$ for some $i\neq j$. The multiple Wiener-It\^o integral
 for a simple function $f$ in the form (\ref{simple}) with respect to the compensated Poisson measure $\widehat{N}$ is defined by
$$I_n(f)=\sum_{k=1}^m\lambda_k\widehat{N}(A_1)\ldots\widehat{N}(A_n).$$
Since the class $\mathcal{S}_n$ is dense in $L^2_{\rm sym}(\mu^n)$ then for every $f\in L^2_{\rm sym}(\mu^n)$, there exists a sequence $\{f_l\}_{l\ge0}\subset\mathcal{S}_n $ such that $f_l\to f$ in $L^2_{\rm sym}(\mu^n)$. Moreover, one can show that $\mathbf{E}[I_n(f_l)^2] = k!\|f_l\|^2$.
Hence, the the multiple Wiener-It\^o integral $I_n(f)$ for a symmetric function $f\in L^2_{\rm sym}(\mu^n)$ can be defined as the limit of the sequence $\{I_n(f_l)\}_{l\ge 0}$ in $L^2(P_N)$ and we denote it as (\ref{integral}).
\begin{proposition}
For $n,m\geq 1$ and $f\in L_{\rm sym}^2(\mu^n)$, $g\in L_{\rm sym}^2(\mu^m)$, then
\begin{enumerate}
  \item $\mathbf{E}[I_n(f)]=0$,
  \item $\mathbf{E}[I_n(f) I_m(g)]=  \delta_{n,m} n!\langle f,g  \rangle_{L^2(\mu^n)}$
\end{enumerate}where $\delta_{n,m}$ is the Kronecker delta.
\end{proposition}

For a measurable function $F:\mathcal{M}(E)\to\overline{\mathbb R}$ and $z\in E$ we define the difference operator as $$D_zF(\eta)=F(\eta+\delta_z)-F(\eta),$$
where $\delta_z$ is the Dirac measure at the point $z$. The iterated difference operator is given by
$$D_{z_1,\hdots,z_n}F=D_{z_1}D_{z_2,\hdots,z_n}F.$$
We define the kernels of $F$  as functions $f_n: E^n\to \overline{\mathbb R}$ given by
$$f_n(z_1,\hdots,z_n)=\frac{1}{n!}\mathbf{E}[D_{z_1,\hdots,z_n}F], n \geq 1, $$
Note that $f_n$ is a symmetric function.\\
We define the Ornstein-Uhlenbeck generator as
$$
LF(\eta) =
 \int\limits_E (F( \eta - \delta_z) - F(\eta)) \eta(dz) - \int\limits_E (F(\eta) - F(\eta + \delta _z))\,  \mu(dz).
$$

\begin{proposition}
For each $F\in L^2(P_N)$, then the kernels $f_n$ are elements of $L^2(\mu^n)$, $n\ge 1$ and uniquely admit the Wiener-It\^o chaos expansion in the form
$$F=\mathbf{E}[F]+\sum_{n=1}^{\infty}I_n(f_n),$$
where the sum converges in $L^2(P_N)$.  Furthermore, for $F,G\in L^2(P_N)$
$${\rm Cov}(F,G)=\mathbf{E}[FG]- \mathbf{E}[F]\mathbf{E}[G]=\sum_{n=1}^\infty n!\langle f_n, g_n \rangle_{L^2(\mu^n)}.$$
\end{proposition}

\begin{proposition}
Let $F\in L^2(P_N)$, and assume that
$$
\sum_{n=1}^\infty n \, n!\|f_n\|^2<\infty.
$$
Then the difference of $F$ at $z\in E$ is given by
$$D_zF=\sum_{n=1}^\infty n I_{n-1}(f_n(y,\cdot)).$$
\end{proposition}

\begin{proposition}
For each random variable $F\in L^2(P_N)$ such that
$$\sum_{n=1}^\infty i^2i!\|f_n\|^2<\infty,$$
then the Ornstein-Uhlenbeck generator $L$ is calculated as
$$LF=-\sum_{n=1}^{\infty}n I_n(f_n).$$
Moreover, its inverse operator is calculated as
$$L^{-1}F=-\sum_{n=1}^\infty\frac{1}{n}I_n(f_n).$$
for each $F\in L^2(P_N)$ such that $\mathbf{E}[F]=0$.
\end{proposition}
For more details of the Malliavin calculus on Poisson space, we refer the reader to \cite{Nualart1990}.

\section{Multi-dimensional Malliavin-Stein inequality}
In the next sequence, we use the probabilistic distance of two $d$-dimensional random vectors $X,Y$  such that $\mathbf{E}(\|X\|_{\mathbb{R}^d}), \mathbf{E}(\|Y\|_{\mathbb{R}^d})<\infty$, which is defined by
$$\Delta(X,Y)=\sup_{g\in \mathcal H}|\mathbf{E}(g(X))-\mathbf{E}(g(Y))|,$$
where $ \mathcal{H}$ is the family of all real-valued functions $g\in C^2(\mathbb{R}^d)$ such that
$$\|g\|_{{\rm Lip}}=\sup_{x\neq y}\frac{|g(x)-g(y)|}{\|x-y\|_{\mathbb{R}^d}}\le 1, \sup_{x\in\mathbb{R}^d}\| {\rm Hess}( g(x))\|\le 1.$$
In the above inequality, $${\rm Hess}(g(z))=\left.\left(\frac{\partial^2 g}{\partial x_i\partial x_j}(z)\right)_{i,j=1}^d\right.$$ stands for the Hessian matrix of $g$ evaluated at a point $z$ and we use the notation of operator norm for a $d\times d$ real matrix $A$  given by $$\|A \| = \sup_{\|x\|_{\mathbb{R}^d}=1} \|Ax\|_{\mathbb{R}^d}.$$

\begin{theorem}[Multi-dimensional Malliavin-Stein inequality]  Consider a random vector  $F=(F_1,\ldots,F_d)\subset L^2(P_N), d\ge2$ such that for $1\leq i\leq d$, $F_i\in
{\rm dom}(D)$, and $\mathbf{E}(F_i)=0$. Suppose that $ X\sim
\mathcal{N}_d(0,C) $, where $C=\{C(i,j): i,j= 1,\ldots,d  \}$ is a
$d\times d $ positive definite symmetric matrix. Then,
$$\begin{array}{ll}\displaystyle
\Delta(F,X) \leq \|C^{-1}\| \|C\|^{1/2} \sqrt{\sum_{i,j=1}^{d} \mathbf{E}[(C(i,j) - \langle  DF_i,-DL^{-1}F_j \rangle_{L^2(\mu)} )^2 ] } \\
\displaystyle
+ \cfrac{\sqrt{2\pi}}{8} \|C^{-1}\|^{3/2} \|C\| \int_E \mu(dz)\mathbf{E}\left[\left(\sum_{i=1}^d|D_z F_i | \right)^2 \left(\sum_{i=1}^d|D_z L^{-1} F_i |  \right) \right].
\end{array}$$
\end{theorem}
For the proof, we refer to \cite{Peccati2010b}.\\

Now we consider a $d$-dimensional random vector $F=(F_1,\ldots,F_d)$ $\subset L^2(P_N)$ with the covariance matrix $\Sigma=\{\Sigma(i,j): i,j= 1,\ldots,d  \}$ such that each component $F_i$ has finite Wiener-It\^o chaos expansions with kernels $f_{i}^{(n)}$, which vanishes if $n>k$. Let give the centered random vector $$G=\sqrt{C\Sigma^{-1}}\left(F-\mathbf{E}[F]\right),$$
where $\sqrt{A}$ stands for the square root of a positive definite matrix $A$, i.e if $A$ has the eigenvalues decomposition $A= P^{-1}{\rm diag}({\lambda_1},\ldots,{\lambda_d})P$
 then$$\sqrt{A}=P^{-1}{\rm diag}(\sqrt{\lambda_1},\ldots,\sqrt{\lambda_d})P.$$
Let us use vector notations $$\nabla_zF=(D_zF_1,\ldots,D_zF_d), \ \nabla_z L F=(D_z LF_1,\ldots ,D_z LF_d)$$ and note that the inequality
\begin{equation}\label{trace}{\rm trace}(AB)\le {\rm trace}(A) {\rm trace}(B)\end{equation}
holds for all positive definite matrices $A,B$.\\
Therefore, by the properties matrix trace ${\rm trace}(AB)={\rm trace}(BA)$ and using the inequality (\ref{trace}), we have
$$\sum_{i,j=1}^{d}(C(i,j) - \langle  DG_i,-DL^{-1}G_j \rangle_{L^2(\mu)} )^2 ={\rm trace}\left[\left(C- \int_E\mu(dz) \nabla_z G(-\nabla_z L G)^T\right)^2\right]$$
$$={\rm trace}\left[\left(\sqrt{C\Sigma^{-1}}\left( \Sigma - \int_E\mu(dz) \nabla_z [F-\mathbf{E}(F)](-\nabla_z L [F-\mathbf{E}(F)])^T\right)\sqrt{C\Sigma^{-1}}^T\right)^2\right]$$
$$={\rm trace}\left[ (C\Sigma^{-1})^2 \left( \Sigma - \int_E\mu(dz) \nabla_z [F-\mathbf{E}(F)](-\nabla_z L [F-\mathbf{E}(F)])^T\right)^2\right]$$
$$\le {\rm trace} [(C\Sigma^{-1})^2]{\rm trace} \left[\left( \Sigma - \int_E\mu(dz) \nabla_z [F-\mathbf{E}(F)](-\nabla_z L [F-\mathbf{E}(F)])^T\right)^2\right]$$
$$=\|C\Sigma^{-1}\|_{F}^2 \sum_{i,j=1}^{d}\left(\Sigma(i,j) - \langle  D[F_i-\mathbf{E}(F)],-DL^{-1}[F_j-\mathbf{E}(F)] \rangle_{L^2(\mu)} \right)^2,$$
where $$\|A\|_F=\sqrt{{\rm trace}(A^TA)}$$ denotes the Frobenius norm of matrix $A$.\\
Note that
$$\Sigma(i,j)={\rm Cov}(F_i,F_j)=\sum_{n=1}^{k}n!\left\langle f_i^{(n)}, f_j^{(n)} \right\rangle_{L^2(\mu^n)},$$
and
$$ \langle  D[F_i-\mathbf{E}(F_i)],-DL^{-1}[F_j-\mathbf{E}(F_j)] \rangle_{L^2(\mu)}=\left \langle \sum_{n=1}^knI_{n-1}(f_i^{(n)}(z,\cdot)),\sum_{n=1}^k I_{n-1}(f_j^{(n)}(z,\cdot))\right\rangle_{L^2(\mu)}.$$
Hence,
$$\begin{array}{l} \displaystyle
\mathbf{E}[(\Sigma(i,j) - \langle  D[F_i-\mathbf{E}(F_i)],-DL^{-1}[F_j-\mathbf{E}(F_j)] \rangle_{L^2(\mu)} )^2]\\ \displaystyle
\le k^2\left(  \sum_{n=1}^{k}\mathbf{E}\left[\left(n!\left\langle f_i^{(n)}, f_j^{(n)} \right\rangle_{L^2(\mu^n)}-n\left \langle I_{n-1}(f_i^{(n)}(z,\cdot)),I_{n-1}(f_j^{(n)}(z,\cdot)) \right\rangle_{L^2(\mu)}\right)^2 \right]\right.\\  \displaystyle
+\left. \sum_{n,m=1,\, n \neq m}^k  \mathbf{E} \left[n^2 \left\langle I_{n-1}(f_i^{(n)}(z,\cdot)),I_{m-1}(f_j^{(m)}(z,\cdot))\right\rangle_{L^2(\mu)}^2\right]\right)\\
 \displaystyle =\ k^2 \sum_{1\le n,m\le k}n^2{\rm Var}\left( \left \langle I_{n-1}(f_i^{(n)}(z,\cdot)),I_{m-1}(f_j^{(m)}(z,\cdot)) \right\rangle_{L^2(\mu)} \right).\end{array}$$
It follows that
$$\sqrt{\sum_{i,j=1}^{d} \mathbf{E}[(C(i,j) - \langle  DG_i,-DL^{-1}G_j \rangle_{L^2(\mu)} )^2 ] }$$
\begin{equation}\label{ineq1}\le k^2 \|C\Sigma^{-1}\|_{F} \sqrt{\sum_{i,j=1}^{d}\sum_{n,m=1}^{k} {\rm Var}\left( \left \langle I_{n-1}(f_i^{(n)}(z,\cdot)),I_{m-1}(f_j^{(m)}(z,\cdot)) \right \rangle_{L^2(\mu)} \right) }.\end{equation}
Morever, by using Holder inequality and the property of matrix norm, we have
$$ \int_E \mu(dz)\mathbf{E}\left[\left(\sum_{i=1}^d|D_z G_i | \right)^2 \left(\sum_{i=1}^d|D_z L^{-1} G_i |  \right) \right]
\le d^{3/2}\int_E \mu(dz)\mathbf{E}\left[\|\nabla_z G \|_{\mathbb{R}^d}^2\|\nabla_zL^{-1} G \|_{\mathbb{R}^d} \right]$$
$$=d^{3/2}\int_E \mu(dz)\mathbf{E}\left[\| \sqrt{C\Sigma^{-1}} \nabla_z [F-\mathbf{E}(F)]\|_{\mathbb{R}^d}^2\|\sqrt{C\Sigma^{-1}} \nabla_zL^{-1}[F-\mathbf{E}(F)]  \|_{\mathbb{R}^d} \right]$$
$$\le d^{3/2}\|\sqrt{C\Sigma^{-1}}\|^3\int_E \mu(dz)\mathbf{E}\left[\|  \nabla_z [F-\mathbf{E}(F)]\|_{\mathbb{R}^d}^2\| \nabla_zL^{-1} [F-\mathbf{E}(F)] \|_{\mathbb{R}^d} \right]$$
$$\le d^{3/2}\|\sqrt{C\Sigma^{-1}}\|^3 \left( \int_E \mu(dz) \mathbf{E}\left[\left(\sum_{i=1}^d|D_z [F_i-\mathbf{E}(F_i)] |^2\right)^2\right]  \right)^{1/2}$$
$$\times \left( \int_E \mu(dz) \mathbf{E}\left[ \sum_{i=1}^d|D_z L^{-1} [F_i-\mathbf{E}(F_i)] |^2   \right]  \right)^{1/2}
$$
$$\le d^2 \|\sqrt{C\Sigma^{-1}}\|^3  \left( \sum_{i=1}^d \int_E \mu(dz) E\left[|D_z [F_i-\mathbf{E}(F_i)] |^4\right]  \right)^{1/2}$$
$$\times \left( \sum_{i=1}^d \int_E \mu(dz) \mathbf{E}\left[|D_z L^{-1} [F_i-\mathbf{E}(F_i)] |^2\right]  \right)^{1/2}$$
$$\le  d^2 \|\sqrt{C\Sigma^{-1}}\|^3  \left( \sum_{i=1}^d   \int\limits_E\mu(dz) k^3 \sum_{n=1}^k  n^4 \mathbf{E}[I_{n-1}(f_i^{(n)}(z,\cdot))^4] \right)^{1/2} $$
$$\times \left( \sum_{i=1}^d                 \int\limits_E\mu(dz)\sum_{n=1}^k \mathbf{E}[I_{n-1}(f_i^{(n)}(z,\cdot))^2]\right)^{1/2} $$
$$=  d^2  \|\sqrt{C\Sigma^{-1}}\|^3 \left( \sum_{i=1}^d  k^3 \sum_{n=1}^k n^4 \mathbf{E}[\|I_{n-1}(f_i^{(n)}(z,\cdot))^2 \|^2] \right)^{1/2}  \left( \sum_{i=1}^d \sum_{n=1}^k  (n-1)!\|f_i^{(n)}\|^2  \right)^{1/2}$$
\begin{equation}\label{ineq2}\le d^{2} k^{7/2} \|\sqrt{C\Sigma^{-1}}\|^3({\rm trace}(\Sigma))^{1/2} \sqrt{\sum_{i=1}^d\sum_{n=1}^k \mathbf{E}\left[\|I_{n-1}(f_i^{(n)}(z,\cdot))^2 \|^2\right]}.\end{equation}
Substituting (\ref{ineq1}) and (\ref{ineq2}) to the inequality in Theorem 3.1 for $G$, we obtain that

\begin{theorem} Let give a $d$-dimensional Gaussian random variable $ X\sim \mathcal{N}_d(0,C)$. Assume that $F=(F_1,\ldots,F_d)$ $\subset L^2(P_N)$ such that ${\rm Cov}(F_i,F_j)=\Sigma(i,j), \ i,j=1,d$ and $F_i$ has finite Wiener-It\^o chaos expansions with kernels $f_{i}^{(n)}$, which vanishes if $n>k$. Then
$$\begin{array}{lll}\displaystyle
\Delta\left(\sqrt{C\Sigma^{-1}}\left(F-\mathbf{E}(F)\right),X\right) \leq \\ \displaystyle
  \cfrac{ \sqrt{2\pi}}{8}d^{2} k^{7/2} \|\sqrt{C\Sigma^{-1}}\|^3 \|C^{-1}\|^{3/2} \|C\| ({\rm trace}(\Sigma))^{1/2} \sqrt{\sum_{i=1}^d\sum_{n=1}^k \mathbf{E}\left[\|I_{n-1}(f_i^{(n)}(z,\cdot))^2 \|^2\right]}\\
\displaystyle +  k^2\|C\Sigma^{-1}\|_{F} \|C^{-1}\| \|C\|^{1/2}\sqrt{\sum_{i,j=1}^{d}\sum_{n,m=1}^{k} {\rm Var}\left( \left \langle I_{n-1}(f_i^{(n)}(z,\cdot)),I_{n-1}(f_j^{(m)}(z,\cdot)) \right \rangle_{L^2(\mu)} \right) }.
\end{array}$$
\end{theorem}

\section{Application for multi-dimensional U-statistics}

In this section we consider the $d$-dimensional vector of U-statistics of the Poisson point process $N$
\begin{equation}\label{ustat}
 F = \left(\sum_{({z}_1,\ldots,{z}_{k_1}) \in S_{k_1}(N)} \phi_1({z}_1,\dots,{z}_{k_1}),\ldots,\sum_{({z}_1,\dots,{z}_{k_d}) \in S_{k_d}(N)} \phi_d({z}_1,\dots,{z}_{k_d}) \right),\end{equation}
where $\phi_i\in L^1_{\rm sym}(\mu^{k_i})$, and $S_{k_i}(N)$ denotes the set of all $k_i$-tuples of distinct points of $N$. This means that each component
$$F_i=\sum_{({z}_1,\dots,{z}_{k_i})\in S_{k_i}(N)} \phi_i({z}_1,\dots,{z}_{k_i})$$
is an  U-statistic of order $k_i$ with respect to the Poisson point process $N$, $i=1,d$.\\

The following properties of (one-dimensional) U-statistics are obtained by Reitzner and Schulte in \cite{Reitzner2011}

\begin{proposition}

Let $F\in L^2(P_N)$ be a U-statistic of order $k$ in the form
$$F=\sum_{({z}_1,\dots,{z}_{k})\in S_k(N)} \phi({z}_1,\dots,{z}_k)$$
Then the kernels of the Wiener-It\^o chaos expansion of $F$ have the form
$$
f_n(z_1,\hdots,z_n)=
\begin{cases}\displaystyle
\binom{k}{n}\int\limits_{E^{k-n}}\phi(z_1,\hdots,z_n,x_1,\hdots,x_{k-n})\, \mu^{k-n}(dx_1,\dots,dx_{k-n}), &n\leq k\\ 0, & n>k.
\end{cases}
$$

\end{proposition}
\begin{proposition}
Assume $F\in L^2(P_N)$, then
\begin{enumerate}
\item If $F$ is a U-statistic, then $F$ has a finite Wiener-It\^o chaos expansion with kernels $f_n\in L^1(\mu^n)\cap L^2(\mu^n) $, $n=1,\hdots,k$.
\item If $F$ has a finite Wiener-It\^o chaos expansion with kernels $f_n\in L^1(\mu^n) \cap L^ 2 (\mu^n)$, $n=1,\hdots,k$, then $F$ is a finite sum of U-statistics and a constant.
\end{enumerate}
\end{proposition}

\begin{proposition} Let $f_i\in \mathcal{S}_{k_i}, \ i=1,\hdots,m$ and $ \Pi$ be the set of all partitions of $Z=\{z_1^{(1)},\dots,z_{n_1}^{(1)},\dots,z_{1}^{(m)} ,\dots, z_{n_m}^{(m)}\}, n_i\le k_i$  such that for each $\pi \in \Pi$,
\begin{enumerate}
 \item $z^{(i)}_{l}, z^{(i)}_{h}\in Z$, $l\neq h$  are always in different subsets of $\pi$, and such that
 \item every subset of $\pi$ has at least two elements.
  \end{enumerate}
 For every partition $\pi\in\Pi$ we define an operator $R^{\pi}$ that replaces all elements of $Z$ in $\prod_{i=1}^m f_i(z_1^{(i)},\hdots,z_{n_i}^{(i)})$ that belong to the same subset of $\pi$ by a new variable $x_j$, $j =1, \hdots,{|\pi|}$, where $|\pi|$ denotes the number of subsets of the partition $\pi$. Then
$$\mathbf{E} \left[\prod_{i=1}^m I_{n_i}(f_i)\right]=\sum_{\pi\in\Pi}\int\limits_{E^{|\pi|}}R^{\pi}(\prod_{i=1}^m f_i(\cdot))(x_1,\hdots,x_{|\pi|})\, \mu^{|\pi|}(dx_1,\dots,dx_{|\pi|}).$$
\end{proposition}

 Using the Proposition 4.3 and the same technique in \cite{Reitzner2011} (Lemma 4.6), we also obtain that if $F=(F_1,F_2,...,F_d)\subset L^2(P_N)$ is a vector of U-statistics in the form (\ref{ustat}) such that $\phi_i, i=1,d$ are simple functions, then all kernels $f_i^{(n)}$ are also simple functions and
$${\rm Var}\left( \left \langle I_{n-1}(f_i^{(n)}(z,\cdot)),I_{m-1}(f_j^{(m)}(z,\cdot)) \right \rangle_{L^2(\mu)} \right) $$
$$= \mathbf{E}\left[ \left \langle I_{n-1}(f_i^{(n)}(z,\cdot)),I_{m-1}(f_j^{(m)}(z,\cdot)) \right \rangle_{L^2(\mu)}^2\right]- \delta_{n,m}\left((n-1)!\left\langle f_i^{(n)}, f_j^{(m)} \right\rangle_{L^2(\mu^n)}\right)^2$$
$$= \int_{E^2}\mathbf{E}\left[ I_{n-1}(f_i^{(n)}(z,\cdot))^2 I_{m-1}(f_j^{(m)}(y,\cdot))^2 \right] \mu^2(dy,dz)$$
$$-\delta_{n,m}\left((n-1)!\left\langle f_i^{(n)}, f_j^{(m)} \right\rangle_{L^2(\mu^n)}\right)^2$$
$$\le \sum_{\pi\in \overline{\Pi}_{n,m}}\int_{E^{|\pi|}}R^{\pi}\left( \left|f^{(n)}_i(.)f^{(n)}_i(.)f^{(m)}_j(.)f^{(m)}_j(.)\right|\right)(x_1,\ldots, x_{|\pi|})\mu^{|\pi|}(dx_1,\ldots, dx_{|\pi|}),$$
and
$$\mathbf{E}\left[\|I_{n-1}(f_i^{(n)}(z,\cdot))^2 \|^2\right]=\int_{E}\mathbf{E}\left[I_{n-1}(f_i^{(n)}(z,\cdot))^4\right]\mu(dz)$$
$$\le \sum_{\pi\in \overline{\Pi}_{n,m}}\int_{E^{|\pi|}}R^{\pi}\left( \left|f^{(n)}_i(.)f^{(n)}_i(.)f^{(n)}_i(.)f^{(n)}_i(.)\right|\right)(x_1,\ldots, x_{|\pi|})\mu^{|\pi|}(dx_1,\ldots, dx_{|\pi|}),$$
where $ \Pi_{n,n}$ stands for the set of partitions satisfying the conditions in Proposition 4.3 with $Z_{n,m}=\{z_{1}^{(1)}\ldots z_{n-1}^{(1)},z_{1}^{(2)}\ldots z_{n-1}^{(2)},z_{1}^{(3)}\ldots z_{m-1}^{(3)},z_{1}^{(4)}\ldots z_{m-1}^{(4)}\}$ and $\overline{\Pi}_{n,m} \subset \Pi_{n,m}$ denotes the set of all partitions in $\Pi_{n,m}$ of such that for any $\pi \in \overline{\Pi}_{n,m}$ and any decomposition of $\{ 1, 2,3,4\}$ into two disjoint sets $M_1, M_2$ there are $i \in M_1, j \in M_2$  and two variables $z_l^{(i)}, z_h^{(j)}$ which are in the same subset of $\pi$.\\
By the formula of kernels in Proposition 4.1., we note that
$$\sum_{\pi\in \overline{\Pi}}\int_{E^{|\pi|}}R^{\pi}\left( \left|f^{(n)}_i(.)f^{(n)}_i(.)f^{(m)}_j(.)f^{(m)}_j(.)\right|\right)(x_1,\ldots, x_{|\pi|})\mu^{|\pi|}(dx_1,\ldots, dx_{|\pi|})$$
$$
\le M_{n,m}(i,j) = \mathbf{1}_{n\le k_i, m\le k_j} \binom{k_i}{n}^2\binom{k_j}{m}^2  \sum_{\pi \in \overline{\Pi}_{n,m}}
\int\limits_{E^{|\pi|}} \int\limits_{E^{2(k_i-m)}}   \int\limits_{E^{2(k_j-n)}}
$$
$$R^{\pi} \left(  |
\prod_{l=1}^2 \phi_i(\cdot , x^{(l)}_1,\hdots,x^{(l)}_{k_i-n})
\prod_{l=3}^4 \phi_j(\cdot , x^{(l)}_1,\hdots,x^{(l)}_{k_j-m})
| \right)(y_1,\hdots,y_{|\pi|}),
$$
\begin{equation}\label{err}
\mu^{|\pi|+2(k_i+k_j-j-i)}(dx_1^{(1)},\hdots,dx_{k_i-n}^{(2)},dx_{1}^{(3)} \ldots dx_{k_j-m}^{(4)}, dy_1, \dots dy_{|\pi|}).\end{equation}
This fact follows that
\begin{theorem}
Assume that $F=(F_1,F_2,...,F_d)\subset L^2(P_N)$ is a vector of U-statistics in the form (\ref{ustat}) such that $\phi_i, i=1,d$ are simple functions. Then
$$\begin{array}{lll}\displaystyle
\Delta\left(\sqrt{C\Sigma^{-1}}\left(F-\mathbf{E}(F)\right),X\right) \leq \\ \displaystyle
  \cfrac{ \sqrt{2\pi}}{8}d^{2} k^{7/2} \|\sqrt{C\Sigma^{-1}}\|^3 \|C^{-1}\|^{3/2} \|C\| ({\rm trace}(\Sigma))^{1/2} \sqrt{\sum_{i=1}^d\sum_{n=1}^k M_{n,n}(i,i)}\\
\displaystyle +  k^2\|C\Sigma^{-1}\|_{F} \|C^{-1}\| \|C\|^{1/2}\sqrt{\sum_{i,j=1}^{d}\sum_{n,m=1}^{k} M_{n,m}(i,j) },
\end{array}
$$
where $k=\max\{k_i, 1\le i\le,d\}$ and $M_{n,m}(i,j), 1\le i,j\le d, 1\le n,m\le k$ are defined in (\ref{err}).
\end{theorem}

Now, we consider that $F=(F_1,F_2,....,F_d)\subset L^2(P_N)$ a vector of $U$-statistics in the form (\ref{ustat}) such that
$$\sum_{({z}_1,\dots,{z}_{k_i})\in S_{k_i}(N)} |\phi_i({z}_1,\dots,{z}_{k_i})|\in L^2(P_N).$$
Then, for each $i=1,2,\ldots,d$ there exists a sequence $\{\phi_{i,l}\}_{l\ge 0}\subset \mathcal{S}_{k_i}$ such that $|\phi_{i,l}|\le |\phi_i|$ and $\phi_{i,l}$ converges to $\phi_i$ $\mu^{k_i}$-almost everywhere. Let give the vector of U-statistics $F^{(l)}=(F_{1,l},\ldots,F_{d,l})$, where
$$F_{i,l}=\sum_{({z}_1,\dots,{z}_{k_i})\in S_{k_i}(N)} \phi_{i,l}({z}_1,\dots,{z}_{k_i}).$$
Hence,
$$|F_{i,l}|\le \sum_{({z}_1,\dots,{z}_{k_i})\in S_{k_i}(N)} |\phi_{i,l}({z}_1,\dots,{z}_{k_i})|\le \sum_{({z}_1,\dots,{z}_{k_i})\in S_{k_i}(N)} |\phi_i({z}_1,\dots,{z}_{k_i})|\in L^2(P_N).$$
Its follow that $F_{i,l}\in L^2(P_N)$, $F_{i,l}$ converges to $F_i$ almost surely and all kernels $f_{i,l}^{(n)}$ in the Wiener-It\^o chaos expansion of $F_{i,l}$ are simple functions.
Note that
$$\Sigma(i,j)={\rm Cov}(F_i,F_j)=$$
$$=\sum_{n=1}^{\infty}n!\binom{k_i}{n}\binom{k_j}{n}\int\limits_{E^n}\ \int\limits_{E^{k_i-n}} \phi_i(z_1,\hdots,z_n,x_1,\hdots,x_{k_i-n})\mu^{k_i-n}(dx_1, \dots dx_{k_i-n})$$ $$\times \int\limits_{E^{k_j-n}}\phi_j(z_1,\hdots,z_n,x_1,\hdots,x_{k_j-n})\, \mu^{k_j-n}(dx_1, \dots dx_{k_j-n}) \mu^n(dz_1, \dots, dz_n).$$
Moreover, the integrals
$$\int\limits_{E^n}\ \int\limits_{E^{k_i-n}} |\phi_i(z_1,\hdots,z_n,x_1,\hdots,x_{k_i-n})|\mu^{k_i-n}(dx_1, \dots dx_{k_i-n})$$
$$\times \int\limits_{E^{k_j-n}}|\phi_j(z_1,\hdots,z_n,x_1,\hdots,x_{k_j-n})|\, \mu^{k_j-n}(dx_1, \dots dx_{k_j-n}) \mu^n(dz_1, \dots, dz_n)$$
always exist for $1\le n\le k_i, 1\le i,j\le d$.
Therefore, by applying the Lebesgue dominated convergence theorem, we obtain that $\Sigma^{(l)}(i,j)\to \Sigma(i,j)$  and $\mathbf{E}(F_{i,l})\to\mathbf{E}(F_{i})$ for $l\to\infty$. Hence,
$$\sqrt{C(\Sigma^{(l)})^{-1}}\left(F^{(l)}-\mathbf{E}(F^{(l)})\right)\to \sqrt{C\Sigma^{-1}}\left(F-\mathbf{E}(F)\right)$$
almost surely for $l\to\infty$.
Note that, the almost sure convergence implies the convergence in the probabilistic distance $\Delta$ and $|M^{(l)}_{n,m}(i,j)|\le |M_{n,m}(i,j)|$ , where $M^{(l)}_{n,m}(i,j)$ is defined when we replace $\phi_{i},\phi_{j}$ by $\phi_{i}^{(l)},\phi_{j}^{(l)}$ in (\ref{err}). Therefore, by using Theorem 4.1 and applying the triangular inequality, we conclude that
\begin{theorem}
Assume that $F=(F_1,\ldots,F_d)\subset L^2(P_N)$ is a vector of U-statistics in the form (\ref{ustat}) such that
$$\sum_{({z}_1,\dots,{z}_{k_i})\in S_{k_i}(N)} |\phi_i({z}_1,\dots,{z}_{k_i})|\in L^2(P_N).$$
Then
$$\begin{array}{lll}\displaystyle
\Delta\left(\sqrt{C\Sigma^{-1}}\left(F-\mathbf{E}(F)\right),X\right) \leq \\ \displaystyle
  \cfrac{ \sqrt{2\pi}}{8}d^{2} k^{7/2} \|\sqrt{C\Sigma^{-1}}\|^3 \|C^{-1}\|^{3/2} \|C\| ({\rm trace}(\Sigma))^{1/2} \sqrt{\sum_{i=1}^d\sum_{n=1}^k M_{n,n}(i,i)}\\
\displaystyle +  k^2\|C\Sigma^{-1}\|_{F} \|C^{-1}\| \|C\|^{1/2}\sqrt{\sum_{i,j=1}^{d}\sum_{n,m=1}^{k} M_{n,m}(i,j) },
\end{array}
$$
where $k=\max\{k_i, 1\le i\le,d\}$ and $M_{n,m}(i,j), 1\le i,j\le d, 1\le n,m\le k$ are defined in (\ref{err}).
\end{theorem}

\begin{corol}
Assume that $\{F^{(l)}\}_{l\ge 0}$ is a sequence of vectors of U-statistics, which are defined as in Theorem 4.2, such that
$$\max_{1\le i,j\le d, 1\le n,m\le k}M_{n,m}^{(l)}(i,j)\to 0$$
for $l\to\infty$, then the law of $\sqrt{C(\Sigma^{(l)})^{-1}}\left(F^{(l)}-\mathbf{E}(F^{(l)})\right)$  converges to the multivariate Gaussian law $ \mathcal{N}_d(0,C)$.
 \end{corol}

\bibliographystyle{amsalpha}
\bibliography{x}
\end{document}